\documentclass[12pt]{amsart}
\usepackage{amscd,amsmath,amssymb,amsfonts}
\theoremstyle{plain}
\newtheorem{thm}{Theorem}
\newtheorem{lem}[thm]{Lemma}
\newtheorem{cor}[thm]{Corollary}
\newtheorem{prop}[thm]{Proposition}
\newtheorem{remark}[thm]{Remark}

\newtheorem{defn}[thm]{Definition}

\numberwithin{thm}{section}
\numberwithin{equation}{section}

\newcommand{\eq}[2]{\begin{equation}\label{#1}#2 \end{equation}}
\newcommand{\ml}[2]{\begin{multline}\label{#1}#2 \end{multline}}

\newcommand{\inj}{\hookrightarrow}

\newcommand{\im}{{\rm im}}

\newcommand{\Tr}{{\rm Tr}}

\newcommand{\sA}{{\mathcal A}}
\newcommand{\sB}{{\mathcal B}}
\newcommand{\sC}{{\mathcal C}}
\newcommand{\sD}{{\mathcal D}}
\newcommand{\sE}{{\mathcal E}}
\newcommand{\sF}{{\mathcal F}}
\newcommand{\sG}{{\mathcal G}}
\newcommand{\sH}{{\mathcal H}}

\newcommand{\sK}{{\mathcal K}}

\newcommand{\sO}{{\mathcal O}}

\newcommand{\E}{{\mathbb E}}

\renewcommand{\H}{{\mathbb H}}

\renewcommand{\P}{{\mathbb P}}

\newcommand{\R}{{\mathbb R}}

\newcommand{\Z}{{\mathbb Z}}

\begin{document}

\title[Relative characters]{Relative algebraic differential characters} 
\author{Spencer Bloch}
\address{Dept. of Mathematics,
University of Chicago,
Chicago, IL 60637,
USA}
\email{bloch@math.uchicago.edu}

\author{H\'el\`ene Esnault}
\address{Mathematik,
Universit\"at Essen, FB6, Mathematik, 45117 Essen, Germany}
\email{esnault@uni-essen.de}
\date{December 1, 1999}
\begin{abstract} Let $f: X \to S$ be a smooth morphism in
characteristic zero, and let $(E, \nabla_{X/S})$ be a regular 
relative connection. We define a cohomology of relative
differential characters on $X$ which receives classes of $(E,
\nabla_{X/S})$. This is a version in  family of the
corresponding constructions performed in \cite{E}. It says in
particular that the partial vanishing of the trace of the
iterated Atiyah
classes can be made canonical. When applied to a family of
curves and $c_2$, the construction yields a connection on 
$f_*c_2(E) \in {\rm Pic}(S)$. Such a connection has been constructed
analytically by A. Beilinson (\cite{Be}).
We also relate our construction to the trace complex as 
defined in \cite{BS}.
\end{abstract}
\maketitle
\begin{quote}
\end{quote}

\section{Introduction}
Let $f: X \to S$ be a smooth family of curves aver a smooth
base $S$ over a field $k$ of characteristic zero, 
and $\nabla_{X/S}: E \to \Omega^1_{X/S}\otimes E$  be a relative
connection. Then the Atiyah class 
$$At(E)\in H^1(X, \Omega^1_X \otimes
{\sE}nd(E))
$$ 
lies in the image of $H^1(X, f^*\Omega^1_S \otimes
{\sE}nd(E))$, and therefore higher Chern classes in $H^i(X, \Omega^i_X)$
die in $H^i(X, \Omega^i_X/f^*(\Omega^i_S))$. In this situation,
when
$k$ is the field of complex numbers, A. Beilinson (\cite{Be})
uses Deligne-Beilinson 
analytic cohomology to show that the images $f_*c_2(E)$ and
$f_*(c_1(E)^2)$, as
analytic bundles, are endowed with canonical holomorphic
connections. 
We show in this note  that a careful
analysis of splitting
principle developed in \cite{E}
leads to the existence of functorial classes
$c_2((E,\nabla_{X/S}))$ and $c_1((E, \nabla_{X/S}))^2$ in the group
$\H^2(X, \sK_2 \xrightarrow{d\log} \Omega^1_S \otimes
\Omega^1_{X/S})$. In other words, the vanishing of the class
in $H^2(X, \Omega^2_X/f^*(\Omega^2_S))$ is made canonical.
The trace of those classes define then
isomorphism classes of line bundles with an algebraic connection
in $\H^1(S, \sO^*_X \xrightarrow{d\log} \Omega^1_S)$.
We show how this example 
is a particular case of a more general theory of relative
algebraic differential characters (see precise formulation in
in section 2, and in subsequent sections for the construction), 
which allows higher dimensional smooth morphisms and higher
classes as well. 
  
In fact, A. Beilinson does
not consider just vector bundles and Chern classes, but rather $G$-principal
bundles under a reductive group $G$ and arbitrary characteristic classes of
weight 2. In this note we give two constructions of classes. The first is
based on a splitting principle and is valid for vector bundles. The
second uses the Weil-algebra homomorphism, as constructed in
\cite{BK}, to define a universal class. The latter should
apply to more general $G$-bundles.

In  a final note, we relate our construction to the trace
complex as defined in \cite{BS}.

It is a pleasure to work out A. Beilinson's idea in a more
algebraic language, hoping that he'll like this framework. It is
also a pleasure to thank him for his generosity.

\section{Relative Cohomology}
Let $f: X \to S$ be a smooth morphism between smooth varieties
over a field $k$ of characteristic 0. Then, for $i\geq 0$, one
defines the subcomplex of the de Rham complex 
\begin{gather} \label{defnf}
\sF^i= f^*\Omega^{i-\bullet}_S\wedge 
\Omega^{2\bullet}_X [-\bullet] \xrightarrow{\iota} \Omega^{\ge i}_X[i], 
\end{gather}
where $f^*\Omega^{i-p}_S\wedge 
\Omega^{2p}_X = \Omega^{i+p}_X$ for $i-p < 0$,
and $\sF^0=\Omega^{\bullet}_X$ for $i=0$.  (Note $\sF^i$ starts with
$f^*\Omega^i_S$ in degree $0$.) The basic object of study is the complex 
\begin{gather} \label{defnad}
\sA\sD_{X/S}^i := {\rm cone} \Big(\sK_i \oplus \sF^i 
\xrightarrow{d\log \oplus -\iota} \Omega^{\ge i}_X[i]\Big)[-1],
\end{gather}
where  the Zariski sheaf $\sK_i$ is defined to be the
image of Milnor $K$ sheaf $\sK^M_i$ in its value at the generic
point $i_{k(X)*}K_i^M(k(X))$, or equivalently, the kernel
of the residue map $$i_{k(X)*}K_i^M(k(X)) \to \oplus_{x\in X^{(1)}}
i_{x*}K_{i-1}^M(k(x)).$$ 

One introduces the
following 
\begin{defn} \label{adr}
The group of relative algebraic differential characters in degree $i$
is the group
\begin{gather}
AD^i(X/S)= \H^i(X, \sA\sD_{X/S}^i).
\end{gather}
\end{defn}
We first define a product structure, following \cite{Be1}, 
the way we did in \cite{E}, which is
compatible to the natural product
\begin{gather} \label{pomega}
\Omega^{\ge i}_X[i] \times \Omega^{\ge j}_X[j] \to \Omega^{\ge (i+j)}_X[i+j],
\end{gather}
the induced product 
\begin{gather} \label{pf}
\sF^i \times \sF^j \to \sF^{(i+j)},
\end{gather} 
on the sub-complexes, and the $K$-product
\begin{gather} \label{pk}
\sK_i \times \sK_j \to \sK_{(i+j)}.
\end{gather}
\begin{defn} \label{product} Let $\alpha \in \R$. We define
\begin{gather}
\sA\sD_{X/S}^i \times \sA\sD^j_{X/S} \to
\sA\sD^{i+j}_{X/S}\notag 
\end{gather}
by
\begin{gather}
\begin{align}
x\cup_\alpha y & = \{ x, y\} & x\in \sK_i, y\in \sK_j \notag\\
\ & = 0 & x\in \sK_i, y \in \sF^i \notag\\
\ &= (1-\alpha) d\log x \wedge y & x\in \sK_i, y \in
\Omega^{\geq j}_X[j] \notag\\
\ & = 0 & x \in \sF^i, y \in \sK_j \notag \\
\ & =x \wedge y & x\in \sF^i, y \in \sF^j \notag\\
\ & = (-1)^{{\rm deg}x}\alpha x \wedge y & x\in \sF^i, y \in
\Omega^{\geq j}_X[j] \notag\\
\ & =\alpha x \wedge d\log y & x\in \Omega^{\geq i}_X[i], y \in
\sK_j\notag\\ 
\ & = (1-\alpha) x \wedge y & x\in \Omega^{\geq i}_X[i], y \in
\sF^j\notag\\ 
\ & =0 & x\in \Omega^{\geq i}_X[i], y\in \Omega^{\geq j}_X[j] \notag
\end{align}
\end{gather}
\end{defn}
\begin{prop} \label{prop:pro}
\begin{itemize}
\item[a)]
For every $\alpha \in \R$, 
the formulae of definition \ref{product} define a product,
compatible with the products \eqref{pomega}, \eqref{pf}, \eqref{pk},
with the graded commutativity rule
\begin{gather}
x \cup_{\alpha}y= (-1)^{{\rm deg}x \cdot{\rm deg}y}y \cup_{(1-\alpha)}
x. \notag
\end{gather}
\item[b)] For all $\alpha, \beta \in \R$, the products
$\cup_{\alpha} $ and $\cup_{\beta} $ are homotopic, so, in
particular, the product on $AD^i(X/S)$ is commutative. 
\end{itemize}
\end{prop}
\begin{proof} One has to show
\begin{gather} \label{comp}
\delta(x\cup_\alpha y)=\delta(x) \cup_\alpha y + (-1)^{{\rm deg}x} x
\cup_{\alpha} \delta(y). 
\end{gather}
Even if the verification is a bit lengthy, we insert it
here, for sake of completeness. 
We have a complex $X$, sum of two complexes
 $X = X_1 \oplus X_2$, with maps $u_i: X_i \to Y$.
We denote by $\varphi= u_1-u_2: X \to Y$. Then we define
${\rm Cone}(\varphi)= X[1]\oplus Y$, with differential
$\delta(x,y)=( -d(x), \varphi(x) + d(y))$, and $\delta[-1]$ on 
${\rm Cone}(\varphi)[-1]=X\oplus Y[-1]$ given by $\delta[-1]=-\delta$, that is
concretely $\delta[-1](x,y)= (d(x), -\varphi(x) - d(y))$.

Now take local sections $a_i\in \sK_i, f_i\in \sF^i, \omega_i\in
\Omega^{\ge i}[i]$ and similarly for $i$ replaced by $j$.
So here, $u_1$ is $d\log$ whereas $u_2$ is the natural embedding.
In order to simplify the notation, we omit the $\alpha$ from the
notations (but not from the computations).

One has
\begin{gather}
\delta( a_i \cup a_j) = -u_1(a_i\cup a_j)=-d\log a_i \wedge d\log a_j\\
\delta(a_i) \cup a_j= -\alpha d\log a_i \wedge d\log a_j
\notag\\ 
a_i \cup \delta(a_j) = (1-\alpha) d\log a_i \wedge (- d\log a_j)
\notag
\end{gather} 
thus \eqref{comp} is satisfied.

One has
\begin{gather}
a_i \cup f_j = 0\\
\delta(a_i) \cup f_j= -(1-\alpha) d\log a_i \wedge f_j\notag\\
\delta(f_j)= (df_j, u_2(f_j)=f_j) \notag\\  
a_i \cup df_j=0\notag\\
a_i \cup u_2(f_j)= (1-\alpha) d\log a_i \wedge f_j \notag
\end{gather}
thus \eqref{comp} is satisfied.

One has
\begin{gather}
\delta(a_i\cup \omega_j) = \delta((1-\alpha)d\log a_i \wedge
\omega_j)) \\= -(1-\alpha)d( d\log a_i \wedge \omega_j)=
(1-\alpha) d\log a_1 \wedge d\omega_j \notag\\
\delta(a_i) \cup \omega_j = 0\notag\\
a_i \cup \delta(\omega_j) = (a_i, -d\omega_j)=(1-\alpha)d\log
a_i \wedge (-d\omega_j) \notag
\end{gather}
thus \eqref{comp} is satisfied.

One has
\begin{gather}
f_i \cup a_j =0 \\
\delta(f_i)= (df_i, u_2(f_i)) \notag\\
df_i \cup a_j= 0 \notag \\
u_2(f_i) \cup a_j = \alpha f_i \cup d\log a_j \notag\\
f_i \cup \delta(a_j) = f_i \cup (-d\log a_j)= (-1)^{{\rm
deg}(f_i)} \alpha f_i \wedge (-d\log a_j)\notag
\end{gather}
thus \eqref{comp} is satisfied.

One has
\begin{gather}
\delta(f_i \cup f_j) = \delta(f_i \wedge f_j)= \\
(df_i \wedge f_j + (-1)^{{\rm deg}(f_i)} f_i \wedge df_j,  u_2(f_i
\wedge f_j)) \in (\sF^{i+j+1}, \Omega^{\ge (i+j)}_X) \notag\\
\delta(f_i) \cup f_j = df_i \cup f_j+  u_2(f_i)\cup f_j =
(df_i \wedge f_j , (1-\alpha) f_i \wedge f_j)\notag\\
f_i \cup \delta(f_j) = f_i \cup df_j + f_i \cup u_2(f_j) =
(f_i \wedge df_j , (-1)^{{\rm deg}(f_i)} \alpha f_i \wedge f_j)
\notag 
\end{gather}  
thus \eqref{comp} is satisfied.

One has
\begin{gather}
\delta(f_i \cup \omega_j)= \delta((-1)^{{\rm deg}(f_i)} \alpha
f_i \wedge \omega_j ) = -d ((-1)^{{\rm deg}(f_i)} \alpha
f_i \wedge \omega_j )\\
\delta(f_i) \cup \omega_j = df_i \cup \omega_j+  u_2(f_i)\cup \omega_j
= (-1)^{{\rm deg}(df_i)} \alpha
df_i \wedge \omega_j +0 \notag\\
f_i \cup \delta(\omega_j)= f_i \cup d\omega_j= 
(-1)^{{\rm deg}(f_i)} \alpha
f_i \wedge (-d\omega_j)\notag
\end{gather}
thus \eqref{comp} is satisfied.

One has
\begin{gather}
\delta(\omega_i \cup a_j) = \delta(\alpha \omega_i \wedge d\log
a_j)= -d(\alpha \omega_i \wedge d\log
a_j)\\
\delta(\omega_i)\cup a_j= -d(\omega_i)\cup a_j = -\alpha
d(\omega_i)\cup d\log a_j \notag\\
\omega_i \cup \delta(a_j)= \omega_i \cup (-d\log a_j) = 0\notag
\end{gather}
thus \eqref{comp} is satisfied.

One has
\begin{gather}
\delta(\omega_i \cup f_j) = \delta((1-\alpha) \omega_i \wedge
f_j))= -d ((1-\alpha) \omega_i \wedge
f_j)) \\
\delta(\omega_i) \cup f_j= -d\omega_i \cup f_j =
-(1-\alpha)d\omega_i \wedge f_j \notag\\
\omega_i \cup \delta(f_j)= \omega_i \cup (df_j, u_2(f_j))\notag
= (1-\alpha) \omega_i \wedge df_j.
\end{gather} 
The degree of $\omega_i$ in ${\rm Cone}(\varphi)[-1]$
equals the degree of $\omega_i$ in the
complex $\Omega^{\ge i}_X$ plus 1. Thus \eqref{comp} is fulfilled.

One has
\begin{gather}
\omega_i \cup \omega_j =0 =\delta(\omega_i)\cup \omega_j =
\omega_i \cup \delta(\omega_j) =0 \notag
\end{gather}
thus  \eqref{comp} is satisfied.
 
The homotopy is defined as follows
\begin{gather}
\begin{align}
h(x\otimes y) & =(-1)^m (\alpha - \beta) x\wedge y & x\in
\Big(\Omega^{\ge i}_X[i]\Big)^{m-1}, y\in \Omega^{\ge j}_X[j]
\notag \\ 
\ & =0 & {\rm otherwise \ }.\notag
\end{align}
\end{gather}
\end{proof}
\section{Splitting Principle}
Let $f: X\to S$ be a smooth proper morphism between smooth
varieties, and let $\nabla_{X/S} : E \to \Omega^1_{X/S} \otimes
E$ be a relative connection on a vectorbundle of rank $r$. 
Let $\pi: \P:=\P(E) \to X$
be the projective bundle associated to $E$. The relative 
connection $\nabla_{X/S}$ defines a splitting
\begin{gather}\label{tau}
\tau : \Omega^1_{\P/S} \to \pi^*\Omega^1_{X/S}
\end{gather}
of the exact sequence
\begin{gather}
0\to \pi^*\Omega^1_{X/S} \xrightarrow{\iota_{/S}} \Omega^1_{\P/S} 
\xrightarrow{p_{/S}} \Omega^1_{\P/X} \to 0,
\end{gather}
such that $\tau \circ \pi^*\nabla_{X/S}$ stabilizes
the tautological sequence
\begin{gather}
0 \to \Omega^1_{\P/X}(1) \to \pi^*E \to \sO_\P(1) \to 0
\end{gather}
(see \cite{E1}, \cite{BE1}, \cite{E}). 
Let $i$ and $p$ be defined by the exact sequence
\begin{gather} \label{ext0}
0\to \pi^*\Omega^1_{X} \xrightarrow{i} \Omega^1_{\P} 
\xrightarrow{p} \Omega^1_{\P/X} \to 0,
\end{gather}
and $q : \Omega^1_{\P} \to \Omega^1_{\P/S}$ be the projection.
\begin{defn}
We define the sheaf
\begin{gather} \label{Omega}
\Omega: = {\rm Ker}\Big( \Omega^1_{\P} \xrightarrow{\tau \circ
q} \pi^*\Omega^1_{X/S}\Big) 
\end{gather}
\end{defn}
{From} the definition, one has an  exact sequence
\begin{gather} \label{ext}
0 \to \pi^*  f^* \Omega^1_S \to \Omega \to \Omega^1_{\P/X}
\to 0
\end{gather} 
\begin{lem}\label{ext=at0}
The extension class 
$$  [\Omega] \in H^1(\P, T_{\P/X}\otimes \pi^*
f^*\Omega^1_S)=H^1(X, {\sE}nd^0(E) \otimes \Omega^1_S)
$$ 
given
by \eqref{ext} is the trace free part of the lifting of 
the Atiyah class of $E$ in $H^1(X, {\sE}nd^0(E)\otimes 
 \Omega^1_X)$ defined by the
choice of $\nabla_{X/S}$. 
\end{lem}
\begin{proof}
The Atiyah sequence can be written
\eq{3.7}{0 \to {\sE}nd(E) \to F \to T^1_X \to 0.
}
Here $F$ is interpreted as infinitesimal symmetries of the bundle $E$. A
connection $\nabla_{X/S}$ relative to $S$ gives rise to an action of
$T^1_{X/S}$, i.e. a lifting $\rho:T^1_{X/S}\to F$ of the inclusion
$T^1_{X/S}\subset T^1_X$. 

The corresponding sequence for infinitesimal
symmetries of the projective bundle $\P = \P(E)$ looks like
$$0 \to {\sE}nd^0(E) \to \pi_* T^1_\P \to T^1_X \to 0
$$
Since symmetries of the vector bundle give rise to symmetries of the
projective bundle, we get a diagram (${\sE}nd^0(E) =
{\sE}nd(E)/k\cdot {\rm id}$) 
$$\begin{CD} @.  @. T^1_{X/S} \\
@. @. @VVV \\
0 @>>> {\sE}nd(E) @>>> F @>>> T^1_X @>>> 0 \\
@. @VVV @VVV @| \\
0 @>>> {\sE}nd^0(E) @>>> \pi_* T^1_\P @>>> T^1_X @>>> 0
\end{CD}
$$
In particular, the bottom line represents the pushout of the Aityah
extension to the tracefree endomorphisms. On the other hand the composite of
the vertical arrows gives rise by adjunction to a map $\pi^*T^1_{X/S} \to
T^1_\P$ which is dual to $\tau$ in \eqref{tau}. It follows that the
reduction of structure of the tracefree Atiyah class defined by the
connection $\nabla_{X/S}$ is given by the bottom line in the diagram
$$\begin{CD}@. @.0 @. 0 \\
@. @. @VVV   @VVV \\
@. @. T^1_{X/S} @= T^1_{X/S} \\
@. @. @VV\rho'V @VVV \\
0 @>>> {\sE}nd^0(E) @>>> \pi_*T^1_\P @>>> T^1_X @>>> 0 \\
@. @| @VVV @VVV \\
0 @>>> {\sE}nd^0(E) @>>> \pi_*\Omega^\vee @>>> f^*T^1_S @>>> 0 \\
@. @. @VVV @VVV \\
@. @. 0 @.  0 
\end{CD}
$$
Given the relation between the map labelled $\rho'$ above and the map
$\tau$ in \eqref{tau}, it is straightforward to check that this bottom line
is obtained from \eqref{ext} by dualizing and pushing forward.  
\end{proof}
We now construct a $\tau$ version of the $\sA\sD$ complex. Certainly, one has
\begin{gather} \label{der}
d (\wedge ^i \Omega) \subset \wedge ^{i-1} \Omega \wedge
\Omega^2_{\P}. 
\end{gather}
\begin{defn} 
Formally replacing $f^*\Omega^i_S$
by $\wedge ^i \Omega$ and $\Omega^{\ge i}_X$ by $\Omega^{\ge
i}_{\P}$ in the definition \ref{defnf}, we define the subcomplex 
\begin{gather}
\sF^i_\tau= \wedge^{i-\bullet}\Omega\wedge 
\Omega^{2\bullet}_{\P} [-\bullet] \xrightarrow{\iota} \Omega^{\ge i}_{\P}[i] 
\end{gather}
of the de Rham complex. Note $\sF^0_\tau=\Omega^{\bullet}_{\P}$ and, in
degree $0$, 
$(\sF^i_\tau)^0 = \wedge^i\Omega$.
\end{defn}
This allows one to define
\begin{defn}
\begin{gather} \label{defnad}   
\sA\sD_{\tau}^i = {\rm cone} \Big(\sK_i \oplus \sF^i_\tau 
\xrightarrow{d\log \oplus -\iota} \Omega^{\ge i}_{\P}[i]\Big)[-1],\notag
\end{gather}
and
\begin{gather}
AD^i_{\tau}(\P)= \H^i(\P, \sA\sD_{\tau}^i). \notag
\end{gather}
\end{defn}
Then one defines a product formally as in definitiom  \ref{product},
replacing $\sA\sD_{X/S}^i$ by $\sA\sD^i_\tau$. To be precise:
\begin{defn} \label{product2} Let $\alpha \in \R$. We define
\begin{gather}
\sA\sD_{\tau}^i \times \sA\sD^j_{\tau} \to
\sA\sD^{i+j}_{\tau}\notag  \notag
\end{gather}
by
\begin{gather}
\begin{align}
x\cup_\alpha y & = \{ x, y\} & x\in \sK_i, y\in \sK_j \notag\\
\ & = 0 & x\in \sK_i, y \in \sF^i \notag\\
\ &= (1-\alpha) d\log x \wedge y & x\in \sK_i, y \in
\Omega^{\geq j}_{\P}[j] \notag\\
\ & = 0 & x \in \sF^i_{\tau}, y \in \sK_j \notag \\
\ & =x \wedge y & x\in \sF^i_{\tau}, y \in \sF^j_{\tau} \notag\\
\ & = (-1)^{{\rm deg}x}\alpha x \wedge y & x\in \sF^i_{\tau}, y \in
\Omega^{\geq j}_{\P}[j] \notag\\
\ & =\alpha x \wedge d\log y & x\in \Omega^{\geq i}_{\P}[i], y \in
\sK_j\notag\\ 
\ & = (1-\alpha) x \wedge y & x\in \Omega^{\geq i}_{\P}[i], y \in
\sF^j_{\tau}\notag\\ 
\ & =0 & x\in \Omega^{\geq i}_{\P}[i], y\in \Omega^{\geq j}_{\P}[j] \notag
\end{align}
\end{gather}
\end{defn}
Of course, proposition \ref{prop:pro} holds true as well,
replacing $AD(X/S)$ by $AD_\tau(\P)$.
\begin{defn} \label{xi}
We denote by $\xi$ the class of the induced partial connection
\begin{gather}
(\sO(1)_{\P}, \tau\circ \nabla_{X/S}) \in \notag \\
AD^1_{\tau}(\P)= \H^1(\P, \sK_1 \xrightarrow{ \tau \circ q \circ
d\log} \pi^*\Omega^1_{X/S}), \notag
\end{gather}
and by $[\xi]$ its image in $\H^1(\P, \sF^1_\tau)$ induced by
the connecting morphism of the exact sequence
\begin{gather}
0 \to \sF^1_\tau[-1] \to 
\Big(\sK_1 \xrightarrow{d\log}
\Omega^1_{\P} \xrightarrow{d} \Omega^2_{\P} \to \ldots \Big)
\to \sA\sD^1_\tau \to 0. \notag 
\end{gather}
\end{defn}
\begin{thm} \label{spl}
The product 
\begin{gather}
AD^i_\tau(\P) \cup \xi^{\cup j} \in AD^{i+j}_\tau(\P) \notag
\end{gather}
induces an isomorphism
\begin{gather}
AD^r_\tau(\P) = AD^r(X/S) \oplus AD^{r-1}(X/S) \cup \xi \oplus
\ldots \oplus AD^1(X/S) \cup \xi^{\cup (r-1)}. \notag
\end{gather}
\end{thm}
\begin{proof} 
{From} the exact sequence
\begin{gather}
\ldots \to \H^{n-1}(\P, \Omega^{\ge i}_{\P}[i]) \to \H^n(\P,
\sA\sD^i_\tau) \to H^n(\P, \sK_i) \oplus \H^n(\P, \sF^i_\tau) \to
\ldots,
\end{gather}
the splitting principle on $H^n(\P, \sK_i)$ and 
$\H^{n-1}(\P, \Omega^{\ge i}_{\P}[i])$,  and the compatiblity of
the products, one just has to see that
\begin{gather} \label{spf}
\H^\ell(\P, \sF^i_\tau)= \H^\ell(X, \sF^i) \oplus \H^{\ell -1}
(X, \sF^{i-1})\cup [\xi] \oplus \ldots H^0(X, \sF^{i-\ell}) \cup
[\xi]^{\cup \ell}  
\end{gather}
for $\ell= i, i-1$. 

To this end, we will define a quasi-isomorphism
\begin{equation}\label{qis}\bigoplus_j \sF^{n-j}[-j] \to
\R\pi_*(\sF_\tau^n).
\end{equation}
\begin{defn}
\begin{itemize}
\item[a)]
 We write the term in degree $q$ in $\sF_\tau^n$ as
$\wedge^{n-q}\Omega\cdot\Omega^{2q}_{\P}\subset
\Omega^{n+q}_\P$. 
\item[b)]
The exact sequence 
\eqref{ext} defines a
decreasing filtration on the exterior powers of $\Omega$, with
$fil^j\wedge^p\Omega$ generated by wedges with at least $j$ entries in
$\pi^*f^*\Omega^1_S$.
\end{itemize}
\end{defn}
\begin{lem}$fil^{p-j}\Omega_X^{p+q-j} \stackrel{\cong}{\to}
R^j\pi_*(\wedge^{p}\Omega\cdot\Omega^{q}_{\P})$.
\end{lem}
\begin{proof}[proof of lemma] By projecting $[\xi]$ we find a
canonical element $\theta\in H^0(X, R^1\pi_*(\Omega))$. Taking powers gives
$\theta^j\in H^0(X, R^j\pi_*(\wedge^j\Omega)$. Since $\pi^*f^*\Omega_S^p
\subset \wedge^p\Omega$, multiplication by $\theta^j$ gives a map as in the
statement of the lemma. To see this map is an isomorphism, we argue by
induction on $q$. Consider the diagram 
\begin{tiny}
$$\minCDarrowwidth.2cm\begin{CD}
R^{j-1}\pi_*(\wedge^p\Omega)\otimes\Omega_{X/S}^{q} @>a>>
R^{j}\pi_*(\wedge^{p+1}\Omega\cdot \Omega_\P^{q-1})
@>>> R^{j}\pi_*(\wedge^p\Omega\cdot \Omega_\P^{q})) @>c>>
\makebox[2cm][l]{$R^{j}\pi_*(\wedge^p\Omega)\otimes\Omega_{X/S}^{q}$} \\
@. @A\cong AbA @A\cup \theta^j AA @A\cong AdA \\
0 @>>> fil^{p+1-j}\Omega_X^{p+q-j} @>>> fil^{p-j}\Omega_X^{p+q-j} @>>>
f^*\Omega^{p-j}_S \otimes \Omega_{X/S}^{q} @>>>  0
\end{CD}
$$
\end{tiny}
 where the top row comes from the short exact sequence of sheaves
\begin{gather}
0 \to \wedge^{p+1}\Omega\cdot \Omega^{q-1}_{\P} \to 
\wedge^{p}\Omega\cdot \Omega^{q}_{\P} \to 
\wedge^{p} \Omega \otimes \Omega^q_{X/S} \to 0 .\notag
\end{gather}
Suppose for a moment we know the lemma is true when $q=0$. It follows that
the maps $d$ in the above diagram are isomorphisms in all degrees. It
follows that the maps $c$ are surjective and the maps $a$ are zero. By
induction the maps $b$ are isomorphisms for all $j$, so it follows that the
middle vertical arrow is an isomorphism as desired. 

It remains to prove the lemma when $q=0$.  One proves by descending induction on $j$ that
$$R^k\pi_*(fil^j\wedge^p\Omega)
\cong \begin{cases} 0 & k>p-j \\ f^*\Omega_S^{p-k} & k\le p-j.
\end{cases}
$$
We have $fil^p\wedge^p\Omega = f^*\Omega_S^p$, so the assertion is clear 
when $j=p$. The induction step comes from consideration of the long exact
sequence
\begin{multline*}R^{k-1}\pi_*(\Omega^{p-j}_{\P/X})\otimes f^*\Omega^j_S \to
R^k\pi_*(fil^{j+1}\wedge^{p}\Omega) \to R^k\pi_*(fil^{j}\wedge^{p}\Omega)
\to \\
\to R^{k}\pi_*(\Omega^{p-j}_{\P/X})\otimes f^*\Omega^j_S
\to R^{k+1}\pi_*(fil^{j+1}\wedge^{p}\Omega),
\end{multline*} 
noting that $R^a\pi_*\Omega^b_{\P/X} \cong \sO_X$ when $a=b$ and is zero
otherwise. Details are omitted. 
\end{proof}

To finish the proof of theorem (\ref{spl}), it suffices to note that the
maps in the lemma are compatible with the de Rham differentials, so we get
quasi-isomorphisms in the derived category, which gives \eqref{qis}. The
desired decomposition \eqref{spf} follows by taking cohomology on $X$.
\end{proof}

\begin{defn} \label{cl}
The decomposition of theorem \ref{spl} allows us  to defines classes
$c_i(E, \nabla_{X/S})\in AD^i(X/S)$ in the usual way as the coefficients
of the equation
\begin{gather}
(\xi)^{\cup r} = \sum_{i=0}^{r-1} (-1)^{(i-1)}c_{r-i}(E,
\nabla_{X/S}) \cup   (\xi)^{\cup i}. \notag
\end{gather}
\end{defn}
\begin{thm}
Let 
\begin{gather}
0 \to (E'', \nabla''_{X/S}) \to (E,\nabla_{X/S}) \to (E', \nabla'_{X/S}) \to 0
\notag 
\end{gather}
be an exact sequence of relative connections. Then one has a
Whitney product formula
\begin{gather}\label{whit}
c_n(E, \nabla_{X/S}) = \sum_{i=0}^{n} c_i(E',
\nabla'_{X/S}))\cup c_{n-i}((E'',
\nabla''_{X/S})
\end{gather} 
\end{thm}
\begin{proof} Let $\P' = \P(E')\stackrel{i}{\inj} \P=\P(E)$, and let $U:=
\P-\P'\stackrel{j}{\inj} \P$. Let $\pi':\P' \to X$ and $p: U \to \P'' =
\P(E'')$ be the natural maps. 
\begin{lem} Define $AD_\tau(U)$ by restricting the cone of complexes of
sheaves used to define $AD_\tau(\P)$ to $U$. One has a pullback map
$$p^* : AD_{\tau''}(\P'') \to AD_\tau(U),
$$
and $p^*(\xi'')=j^*(\xi)$. Finally, in degree $r''=\text{rank}(E'')$ we have
$$\ker \Big(j^*:AD_\tau^{r''}(\P) \to AD_\tau^{r''}(U)\Big)\cong \Z.$$ 
\end{lem}
\begin{proof}[proof of lemma] Compatibility of the connections on $E$ and
$E''$ leads to a commutative diagram of sheaves on $U$
$$\begin{CD}@. 0  @. 0 \\
@. @VVV @VVV \\
0 @>>> p^*\Omega_{\tau''} @>>> p^*
\Omega^1_{\P''} @>>> p^*\pi''{}^*\Omega^1_{X/S} @>>> 0  \\
@. @VVaV @VVa\makebox[.2cm][l]{$\qquad       (*)$}V @| \\
0 @>>> \Omega_\tau|_U @>>> \Omega^1_U @>>> p^*\pi''{}^*\Omega^1_{X/S} @>>> 0
\\
@. @VVV @VVV \\
@. \Omega^1_{U/\P''} @= \Omega^1_{U/\P''} \\
@. @VVV @VVV \\
@. 0  @. 0 
\end{CD}
$$
Here of course we use the notation $\Omega_\tau$,
$\Omega_{\tau ''}$ etc to indicate in which projective bundle we
consider the construction \ref{Omega} of $\Omega$.
Exterior powers of the maps labelled $a$ enable one to define
$$p^*:AD_{\tau''}(\P'') \to AD_\tau(U).$$ 
Commutativity of the
square labelled 
$(*)$ above implies, using definition (\ref{xi}) that $p^*(\xi'')=j^*(\xi)$. 

Finally, regarding the kernel of $j^*$ in degree $r''$, we show here that it
is a subgroup of $\Z$. Once we describe the Gysin homomorphism below,
it will be clear this kernel is nonzero. Since $\P'\subset \P$ has
codimension
$r''$, purity results for local cohomology imply
$$\H^{r''}_{\P'}(\P, \sA\sD^{r''}) \inj \ker(\H^{r''}_{\P'}(\P,\sK_{r''}\oplus
\wedge^{r''}\Omega_\tau) \to H^{r''}_{\P'}(\P,\Omega^{r''}_\P)).
$$
Further we have
$$H^{r''}_{\P'}(\P,\sK_{r''}) \cong \Z;\quad
H^{r''}_{\P'}(\P,\wedge^{r''}\Omega_\tau)\inj
H^{r''}_{\P'}(\P,\Omega_\P^{r''}).
$$
It follows that the local cohomology of $\sA\sD_\tau$, 
which maps onto the kernel of
$j^*$, is a subgroup of $\Z$ as claimed. 
\end{proof}
\begin{lem} Let $i:\P' \inj \P$ be the inclusion. There is defined a Gysin
map
$$i_* : AD_{\tau'}^k(\P') \to AD_\tau^{k+r''}(\P).
$$
The image $i_*(AD^0_{\tau'}(\P')) = \Z$ and $1$ maps
to the cycle class
of $\P^{r'}$ in $CH^{r ''}(\P)$. There is a projection formula
$$i_*(x\cdot i^*y) = i_*(x)y. 
$$
\end{lem}
\begin{proof}We have a commutative diagram
$$\begin{CD}@. 0 @.  0 \\
@. @VVV @VVV \\
@. N^\vee_{\P'/\P} @= N^\vee_{\P'/\P} \\
@. @VVV @VVV \\
0 @>>> i^*\Omega_\tau @>>> i^*\Omega^1_\P @>>> \pi'{}^*\Omega^1_{X/S} @>>> 0
\\
@. @VVV @VVV @| \\
0 @>>> \Omega_{\tau'} @>>> \Omega^1_{\P'} @>>> \pi'{}^*\Omega^1_{X/S} @>>>
0\\
@. @VVV @VVV \\
@. 0 @. 0
\end{CD}
$$
The left column leads to maps 
$$\det(N^\vee_{\P'/\P})\otimes\wedge^p\Omega_{\tau'} \to
i^*\wedge^{p+r''}\Omega_\tau.
$$
On the other hand, Grothendieck duality theory, \cite{G}, gives
$$\det(N_{\P'/\P})\otimes i^*\wedge^{p+r''}\Omega_\tau \cong \sE
xt^{r''}_{\sO_\P} (\sO_{\P'},\wedge^{p+r''}\Omega_\tau) \to
\sH^{r''}_{\P'}(\P,\wedge^{p+r''}\Omega_\tau), 
$$
where $\sH_{\P'}$ is the Zariski local cohomology sheaf. Composing these
arrows, we get
$$\wedge^p\Omega_{\tau'} \to \sH^{r''}_{\P'}(\P,\wedge^{p+r''}\Omega_\tau).
$$
An elaboration on this construction, using the middle column of the previous
commutative diagram and the analogous map on $K$-sheaves
$$\sK_{p,\P'} \to \sH^{r''}_{\P'}(\sK_{p+r'',\P}) 
$$
gives a map of complexes
$$\sA\sD^p_{\P'} \to \sH^{r''}_{\P'}(\sA\sD^{p+r''}).
$$
Finally, using purity, we find 
\begin{multline*}AD^p_{\tau '}(\P') = \H^p(\P',\sA\sD^p_{\tau ' }) \to
H^p(\P',\sH^{r''}_{\P'}(\sA\sD^{p+r''}_{\tau})) \to \\
\H^{p+r''}_{\P'}(\P,\sA\sD^{p+r''}_{\tau })\to 
\H^{p+r''}(\P,\sA\sD^{p+r''}_{\tau})=AD^{p+r''}_{\tau}(\P). 
\end{multline*}
In the special case $p=0$, we find
\begin{equation*}\Z = AD^0(\P') \to \H^{r''}_{\P'}(\P,\sA\sD^{r''}_{\tau}) \to
AD^{r''}_{\tau}(\P). 
\end{equation*}
Composing with the map $AD^{r''}_{\tau }(\P) \to CH^{r''}(\P)$, we see that the
above map is injective, and in fact is an isomorphism,
taking 1 to the class of $\P'$.  

Finally, the projection formula is a straightforward consequence of the
multiplication on the complex $\sA\sD_\tau$, which makes
$AD_{\tau' }(\P')$ an
$AD_\tau(\P)$-module. 
\end{proof}

The proof of the Whitney formula \eqref{whit} is now straightforward. Write
 $$F'(\xi)=\sum (-1)^{r'-i}c_{r'-i}(E',\nabla')\xi^i,\
F''(\xi)=\sum (-1)^{r''-i}
c_{r''-i}(E'',\nabla'')\xi^i$$ for the polynomials associated to $E'$ and
$E''$. We know that 
$$j^*F''(\xi) = F''(\xi'')=0,$$
so necessarily $F''(\xi) = i_*(1)$. Also, $F'(\xi')=0$, so
$$0 = i_*(F'(\xi')) = i_*(1)F'(\xi)=F''(\xi)F'(\xi).
$$
Since $F(\xi) = \sum (-1)^{r'+r''- i}
c_{r'+r''-i}(E,\nabla)\xi^i$ is the unique monic
polynomial in $\xi$ of degree $r'+r''$ vanishing in $AD_\tau(\P)$, we conclude
$F=F'F''$. 
\end{proof}

\begin{cor} 
Let $Y \subset X$ be a
fiber of $f$. The
classes $c_i(E, \nabla_{X/S}) \in AD^i(X/S)$ specialize to
the classes $c_i(E|_{Y},\nabla_{X/S}|_{Y}) \in AD^i(Y) $
defined in \cite{E}.
\end{cor}

Let $d$ be the relative dimension of the morphism $f$. The existence for
$n-a<d$ of a trace or transfer map
$$\R f_*\Big(\Omega^n_X/f^*\Omega^a_S\cdot \Omega^{n-a}_X) \to
\Omega^{n-d}_S[-d]
$$
leads to a transfer map
\begin{gather}
f_*: AD^n(X/S) \to \H^{n-d}(S, \sK_{n-d} \xrightarrow{d\log}
 \Omega^{n-d}_S \xrightarrow{d}
\ldots \to \Omega^{n-[\frac{d}{2}]-1}_S)
\end{gather}
defined for example by taking the Gersten-Quillen resolution of
the $\sK$ sheaves, and the Cousin resolution of the coherent
sheaves of forms. 
In particular
\begin{cor} \label{conn}
\begin{gather} 
f_*(c_{d+1}(E, \nabla_{X/S})) \in \H^1(S, \sK_1
\xrightarrow{d\log} \Omega^1_S \xrightarrow{d} \ldots \to
\Omega^{[\frac{d+1}{2}]}_S)\notag 
\end{gather}
is the isomorphism class of a rank one line bundle with a
connection, which is flat for $d \ge 3$.
 \end{cor}

\section{Universal construction via the Weil algebra}

In this section we give another construction of the relative classes, using
unpublished work of A. Beilinson and D. Kazhdan \cite{BK}. These authors 
define a filtered differential graded algebra
$\Omega^{\bullet}_{X,E}\supset F^n\Omega^{\bullet}_{X,E} $,
such that $(\Omega^{\bullet}_X,
\Omega^{\ge n}_X) \xrightarrow{\cong}  (\Omega^{\bullet}_{X,E},
F^n\Omega^{\bullet}_{X,E})$, together with a Weil homomorphism
$S^n (\sG^*)^G \xrightarrow{w^n(E)} F^n\Omega^{\bullet}_{X,E}[2n]$. For us,
the algebraic group is $G=GL(r)$, and $S^n (\sG)^G$ are the
$G$-invariant polynomials on the Lie algebra $\sG$ over $k$.  

We begin by recalling the Beilinson-Kazhdan construction. Let $p:\E\to X$ be
the
$GL(r)$-torsor corresponding to a vector bundle
$E$ on
$X$. Define
$$\Omega^1_{X,E} := (p_*\Omega^1_{\E})^{GL(r)}.
$$
The Atiyah sequence can be reinterpreted as the exact sequence of
$GL(r)$-invariant relative differentials pushed down to $X$
\eq{4.1}{0 \to \Omega^1_X \to \Omega^1_{X,E}
\stackrel{\pi}{\longrightarrow}\sE nd(E) \to 0.
}
Given an exact sequenc of vector bundles $0\to A\to B\stackrel{\pi}{\to} C\to
0$, the generalized Koszul sequence yields for any $n$
\ml{4.2}{0 \to \wedge^nA \to \wedge^nB \stackrel{\delta}{\longrightarrow}
\wedge^{n-1}B\otimes C \stackrel{\delta}{\longrightarrow} \\
\wedge^{n-2}B\otimes Sym^2(C)\to\ldots\to Sym^n(C) \to 0 
}
with 
$$\delta(b_1\wedge\ldots\wedge b_p\otimes c_1\cdot\ldots\cdot
c_{n-p}) = \sum_i (-1)^{i-1}b_i\wedge\ldots
\widehat{b_i}\ldots\wedge b_p\otimes \pi(b_i)\cdot c_1\cdot\ldots\cdot
c_{n-p}
$$
Combining \eqref{4.1} and \eqref{4.2} we define a complex
$\Omega^\bullet_{X,E}$ to be the simple complex associated to the first
quadrant double complex: 
\eq{4.3}{\minCDarrowwidth.2cm \begin{CD}@. @. 0 @>>> \ \cdots
\\ @. @. @AAA \\
@. 0 @>>> Sym^2(\sE nd(E)) @>>>
\ \cdots
\\
 @. @AAA @AAA \\
0 @>>> \sE nd(E) @>d'>> \sE nd(E)\otimes \Omega^1_{X,E} @>>>\ \cdots\\
@AAA  @A d''=\pi AA @AAd'' A \\
\sO_X @>d'>>  \Omega^1_{X,E} @>d'>>\wedge^2 \Omega^1_{X,E} @>>>\ \cdots
\end{CD}
}
Grading by total degree gives
\eq{4.4}{\Omega^n_{X,E} = \bigoplus_{\substack{a+b=n\\a\ge
b}}\wedge^{a-b}\Omega^1_{X,E}\otimes Sym^b(\sE nd(E)) := \bigoplus
\Omega_{X,E}^{a,b}  }
The $n$-th column in \eqref{4.3} is a  resolution of $\Omega^n_X$, so
defining 
\eq{4.5}{F^p\Omega^n_{X,E} := \bigoplus_{\substack{a+b=n \\a\ge
p}}\Omega^{a,b}_{X,E} 
}
we get a filtered quasi-isomorphism
\eq{4.6}{(\Omega^\bullet_X, \Omega_X^{\ge p})
\stackrel{\cong}{\to}(\Omega_{X,E}^\bullet, F^p\Omega_{X,E}^\bullet).   
}
In addition, $\sE nd(E)$ is a twisting of the Lie algebra $\frak g$ of
$GL(r)$, so there is a map
\eq{4.7}{w^n(E):Sym^n(\frak g^\vee)^{GL(r)} \to Sym^n(\sE nd(E)) \to
F^n\Omega^\bullet_{X,E}[2n]. 
}

Assume now we are in a relative situation, with $f:X \to S$ and
$(E,\nabla_{X/S})$ as in section 2. 
\begin{lem}The connection $\nabla_{X/S}$ determines a descent of the Atiyah
extension \eqref{4.1} to an extension
\eq{4.8}{ 0 \to f^*(\Omega^1_S) \to \Omega^1_\nabla \to \sE nd(E) \to 0.
}
I.e. the Atiyah extension comes from \eqref{4.8} by pushout $f^*\Omega^1_S
\to \Omega^1_X$. 
\end{lem}
\begin{proof} A connection defines an infinitesimal action of vector fields
over $X$ on $E$, i.e. a splitting of the Atiyah sequence \eqref{3.7}. Thus,
a relative connection leads to a diagram
\eq{4.9}{\begin{CD}0 @>>> \Omega^1_X @>>> \Omega^1_{X,E} @>>> \sE nd(E) @>>>
0 \\
@. @VVV @VVhV @| \\
0 @>>> \Omega^1_{X/S} @>>> \Omega^1_{X/S,E}
@>\stackrel{\sigma}{\longleftarrow}>>
\sE nd(E) @>>> 0 
\end{CD}
}
The descent comes by defining
\eq{4.10}{\Omega^1_\nabla := \{z\in \Omega^1_{X,E}\ | \ h(z) \in
\im(\sigma)\}. 
}
\end{proof}

Recall the complex $\sF^n$ defined in \eqref{defnf}. Our next objective
is to mimick the above Weil construction, replacing \eqref{4.1} by
\eqref{4.8}. The construction leads to complexes $\sF^n_\nabla$ and
quasi-isomorphisms $\sF^n \inj \sF^n_\nabla$ analogous to \eqref{4.6}. Note
by \eqref{4.10} that $\Omega^1_\nabla \subset \Omega^1_{X,E}$, so one may
define $\wedge^a\Omega^1_\nabla\cdot\wedge^b\Omega^1_{X,E} \subset
\wedge^{a+b}\Omega^1_{X,E}$. 
\begin{lem}\label{lem4.2}Recall
$(\sF^n)^j=f^*\Omega^{n-j}_S\cdot\Omega^{2j}_X. 
$
There is a Koszul type resolution
\ml{4.11}{0 \to (\sF^n)^j \to
\wedge^{n-j}\Omega^1_\nabla\cdot\wedge^{2j}\Omega^1_{X,E} \to \\
\big(\wedge^{n-j-1}\Omega^1_\nabla\cdot\wedge^{2j}\Omega^1_{X,E}\big)\otimes
\sE nd(E)
\to \ldots  \\
\to\big( 
\wedge^{n-j-p}\Omega^1_\nabla\cdot\wedge^{2j}\Omega^1_{X,E}\big)\otimes
Sym^p(\sE nd(E)) \to\ldots \\
\to Sym^{n+j}(\sE nd(E)) \to 0 
}
\end{lem}
\begin{proof}We start with the resolution
\ml{4.12}{0 \to \Omega^q_X \to \wedge^q\Omega^1_{X,E} \to
\wedge^{q-1}\Omega^1_{X,E}\otimes \sE nd(E) \to \ldots \\
 \to Sym^q(\sE nd(E)) \to 0  }
which we filter
\begin{gather}fil^i(\Omega^q_X) = f^*\Omega^i_S\cdot\Omega^{q-i}_X; \\
fil^i \big(\wedge^{q-k}\Omega^1_{X,E}\otimes Sym^k(\sE nd(E)) \big) =  \notag
\\ =\begin{cases}
\big(
\wedge^{i-k}\Omega^1_\nabla\cdot\wedge^{q-i}\Omega^1_{X,E}\big)\otimes
Sym^k(\sE nd(E)) & i\ge k \\
0 & i<k\end{cases} \notag
\end{gather}
This filtration is compatible with the differential, and (using
$\Omega^1_{X,E}/\Omega^1_\nabla \cong \Omega^1_{X/S}$), we find
\ml{4.14}{gr^i\big(\wedge^{q-k}\Omega^1_{X,E}\otimes Sym^k(\sE nd(E)) \big)
\\ =\begin{cases}\wedge^{i-k}\Omega^1_\nabla\otimes
\Omega^{q-i}_{X/S}\otimes Sym^k(\sE nd(E)) & i\ge k \\
0 & i<k  \end{cases}
}
The complex $gr^i$ is obtained by tensoring the Koszul resolution
$$0 \to f^*\Omega^i_S \to \wedge^i\Omega^1_\nabla \to 
\wedge^{i-1}\Omega^1_\nabla\otimes \sE nd(E)  \to \ldots 
$$
with $\Omega^q_{X/S}$ and is thus exact. It follows that $fil^i$ is a 
resolution of $f^*\Omega^i_S\cdot\Omega^{q-i}_X$. Taking $i=n-j;\ q= n+j$,
we obtain the desired resolution of $(\sF^n)^j$. 
\end{proof}

\begin{defn} The complex
\eq{4.15}{(\sF^i_\nabla)^\bullet = \bigoplus_{\substack{a+b=\bullet \\a\ge
b}} (\sF^i_\nabla)^{a,b} \subset \bigoplus_{\substack{a+b=\bullet \\a\ge
b}} \Omega_{X,E}^{a,b}
}
is defined by 
\eq{4.16}{(\sF^i_\nabla)^{a,b} = \begin{cases}
\big(\wedge^{2i-a-b}\Omega^1_\nabla\cdot\wedge^{2a-2i}\Omega^1_{X,E}\big)
\otimes Sym^b(\sE nd(E)) & 2i\ge a+b;\ a\ge i \\
0 & i>a \\
\Omega^{a,b}_{X,E} & a\ge i;\ 2i\le a+b \end{cases}.   
}
The differentials are induced by the differentials on
$\Omega^{\bullet}_{X,E}$. 
\end{defn}

\begin{remark}The content of lemma (\ref{lem4.2}) is that one has a
resolution
$$ 0 \to \sF^i[-i] \to (\sF_\nabla^i)^{i,0} \to  (\sF_\nabla^i)^{i,1}\to
\ldots 
$$
with differential $d''$. 
\end{remark}
\begin{remark} The Beilinson-Kazhdan-Weil homomorphism $w_n(E)$, \eqref{4.7},
takes values in $(\sF^n_\nabla)^{n,n}=\Omega^{n,n}_{X,E} = Sym^n(\sE
nd(E))\subset F^n\Omega^{2n}_{X,E}$. 
\end{remark}
In particular, the complex 
\begin{gather}
{\rm cone } \Big( \sK_n \oplus S^n(\sG^*)^G[-n]
\xrightarrow{d\log \oplus w_n(E)}
F^n\Omega^{\bullet}_{X, E}[n]\Big) [-1]
\end{gather}
used in \cite{E} to define the absolute differential
characters is endowed with a map
\begin{gather}
{\rm cone } \Big( \sK_n \oplus S^n(\sG^*)^G[-n]
\xrightarrow{d\log \oplus w_n(E)}
F^n\Omega^{\bullet}_{X, E}[n]\Big) [-1]
\\
\xrightarrow{\iota_\nabla=(1\oplus w_n, 1)}
{\rm cone } \Big( \sK_n \oplus \sF^n_\nabla [n]
\xrightarrow{d\log \oplus w_n(E)}
F^n\Omega^{\bullet}_{X, E}[n]\Big) [-1]\notag
\end{gather}
where the latter complex is quasi-isomorphic to
$\sA\sD^n_{X/S}$.

Let $X =\cup_i X_i$ be a Zariski covering trivializing $E$ and
let $[E]: X_{\bullet} \to BG_{\bullet}$ be the induced map to the
simplicial classifying scheme $BG_{\bullet}$. We denote by $E_{{\rm un}}$ the 
universal bundle on $BG_{\bullet}$.  
Since 
\begin{gather}
\H^n(BG_\bullet, {\rm cone } \Big( \sK_n \oplus S^n(\sG^*)^G[-n]
\xrightarrow{d\log \oplus w_n(E_{{\rm un}})}
F^n\Omega^{\bullet}_{BG_{\bullet}, E_{{\rm un}}}[n]\Big) [-1])\\
= H^n(BG_\bullet, \sK_n) \notag
\end{gather}
(see \cite{E}), one has a well defined universal class
\begin{multline}
c_{n, {\rm an}}\in \\
\H^n(BG_\bullet, {\rm cone } \Big( \sK_n
\oplus S^n(\sG^*)^G[-n] 
\xrightarrow{d\log \oplus w_n(E_{{\rm un}})}
F^n\Omega^{\bullet}_{BG_{\bullet}, E_{{\rm un}}}[n]\Big) [-1]).
\end{multline}
One way to define $c_{n, {\rm an}}$ 
is again via the splitting principle on $BG_\bullet$.
Via the functoriality map
\begin{gather}
[E]^*: {\rm cone } \Big( \sK_n \oplus S^n(\sG^*)^G[-n]
\xrightarrow{d\log \oplus w_n(E_{{\rm un}})}
F^n\Omega^{\bullet}_{BG_{\bullet}, E_{{\rm un}}}[n]\Big) [-1] \notag \\
\to R[E]_*{\rm cone } \Big( \sK_n \oplus S^n(\sG^*)^G[-n]
\xrightarrow{d\log \oplus w_n(E)}
F^n\Omega^{\bullet}_{X_{\bullet}, E}[n]\Big) [-1] \notag
\end{gather}
this defines a class
\begin{multline}
[E]^*(c_{n, {\rm un}})\in \\
\H^n(X_\bullet, {\rm cone } \Big( \sK_n \oplus
S^n(\sG^*)^G[-n]
\xrightarrow{d\log \oplus w_n(E)}
F^n\Omega^{\bullet}_{X_{\bullet}, E}[n]\Big) [-1])
\end{multline}
and via $\iota_\nabla$, one obtains a class
\begin{defn}
\begin{gather}
c'_n(E, \nabla_{X/S}):=\iota_\nabla ([E]^*(c_{n, {\rm un}})) \in
AD^n(X/S).\notag \end{gather}
\end{defn}

It remains to compare those classes to the classes $c_n(E,
\nabla_{X/S})$ constructed in the previous section with the
splitting principle.
As in \cite{E}, theorem 3.26, one has to verify the Whitney 
formula (4). One proceeds as in (4) of loc. cit.
\begin{prop}
One has
\begin{gather}
c_n(E, \nabla_{X/S}) =
c'_n(E, \nabla_{X/S}) \in AD^i(X/S) \notag
\end{gather}
\end{prop}
\section{The image of $c_2$ in a family of curves} 

We now specialize the discussion to the case where $f: X \to S$
is a smooth family of curves. Let $(E, \nabla_{X/S})$ be a
relative connection, to which one has assigned classes
\begin{gather}
c_2(E, \nabla_{X/S}), c_1(E, \nabla_{X/S})^2 \in AD^2(X/S).
\end{gather}
One considers the trace map
\begin{gather} \label{tr}
f_*: AD^2(X/S) \to AD^1(S/k) 
\end{gather}
introduced in corollary \ref{conn}.
Since $AD^1(S/k)$ is the group
of isomorphism classes of line bundles with connection,
we have constructed in this way a connection on the line bundles
\begin{gather}
f_*c_2(E) \ {\rm and \ } f_*(c_1(E)^2). \notag
\end{gather}

We want to identify the curvature of the connection.

We consider the exact sequence of complexes on $X$
\begin{gather} \label{gm}
0 \to \Big(f^*\Omega^2_S \to f^*\Omega^2_S \otimes
\Omega^1_{X/S}\Big)[-1]\\ 
\to \Big(\sK_2 \xrightarrow{d\log} \Omega^2_X \to
\Omega^3_X/f^*\Omega^3_S\Big) \to \sA\sD^1_{X/S} \to 0 .\notag
\end{gather} 
(Recall that $f$ has dimension 1, thus $\Omega^3_X/f^*\Omega^3_S=
f^*\Omega^2_S\otimes \Omega^1_{X/S}$.) Then given a class
$\gamma \in AD^2(X/S)$, the connecting morphism $\delta$ of \ref{gm}
defines a class 
\begin{gather}
\delta(\gamma) \in \H^2(X,  f^*\Omega^2_S \to f^*\Omega^2_S \otimes
\Omega^1_{X/S}) \stackrel{f_*}{\to} \\ H^0(S, \Omega^2_S \otimes
R^2f_*(\Omega^{\bullet}_{X/S})) = H^0(S, \Omega^2_S).\notag 
\end{gather}
\begin{prop}
Let $f_*(\gamma)=
(L, \nabla) \in AD^1(S/k)$, with curvature $\nabla^2 \in H^0(S,
\Omega^2_S)$. Then one has
\begin{gather}
\nabla^2 = f_*\delta(\gamma). \notag
\end{gather}
\end{prop}  
\begin{proof}Briefly, one has diagrams of complexes of sheaves on $X$
\eq{5.5}{\begin{CD}@. \sK_{2,X} @= \sK_{2,X} \\
@. @VVV @VVV \\
f^*\Omega^2_S @>>> \Omega^2_X @>>> \Omega^1_S\otimes\Omega^1_{X/S} \\
@VVV @VVV @. \\
f^*\Omega^2_S\otimes\Omega^1_{X/S} @>\cong >> \Omega^3_X/f^*\Omega^3_S, 
\end{CD}
}
and on $S$
\eq{5.6}{\begin{CD} @. \sO_S^\times @= \sO_S^\times \\
@. @VVV @VVV \\
@. \Omega^1_S @= \Omega^1_S \\
@. @VVV @. \\
\Omega^2_S @= \Omega^2_S. 
\end{CD}
}
Using Gersten resolutions and Cousin complexes, one defines a trace map 
$$\R f_*(\text{diagram \eqref{5.5}})[1] \to \text{diagram \eqref{5.6}},
$$
which yields a commutative diagram
$$\begin{CD}AD^2(X/S) @>\delta>> \H^2(X, f^*\Omega^2_S \to
f^*\Omega^2_S\otimes\Omega^1_{X/S}) \\
@VVf_* V @VV f_* V \\
AD^1(S) @>d>> \Gamma(S,\Omega^2_S).
\end{CD}
$$

\end{proof}

\subsection{Beilinson's gerbe construction}
Now let us explain the construction (see \cite{Be}),
which is the whole motivation for this note: A. Beilinson 
constructs such a connection on $f_{*}c_2(E)$ (and
more generally on the image of classes of $G$-bundles of weight
2) in the analytic Deligne cohomology.
We want to
compare our construction to the one in \cite{Be}.
We should emphasize that this comparison lacks precision
in two points.
First, we give our own interpretation of Beilinson's construction
in terms of a specific cocyle, and then we do not give the
details of how this precise cocyle yields the same class as ours. 

Beilinson uses gerbes as follows. Consider the Atiyah torsor
of $E$. Given a local trivialization of $E$ on $X=\cup_i X_i$,
this defines transition functions $g_{ij} \in \sC^1({\sA}ut(E_{X}))$, and 
$dg_{ij}g^{-1}_{ij} \in \sC^1(\Omega^1_X \otimes {\sE}nd(E_{X}))$
is a cocyle for the Atiyah torsor. Then the $\Omega^2_X$-gerbe 
$c_1^2-2c_2$ is
represented by the cocyle ${\rm Tr}(dg_{ij}g_{ij}^{-1}
dg_{jk}g_{jk}^{-1}) \in \sC^2(\Omega^2_{X})$
whereas the $\Omega^2_X$-gerbe $c_1^2$ is represented by the
cocyle $({\rm Tr}dg_{ij}g_{ij}^{-1})^2 \in \Omega^2_{X_{ijk}}$.
Let us consider the Cech resolution ($j_I : X_I \inj X$)
\begin{gather}
\Omega^2_X \to \prod_i j_{i*}\Omega^2_{X_i} \xrightarrow{d_0}
\prod_{i<j}j_{ij*} \Omega^2_{X_{ij}} \xrightarrow{d_1}
\prod_{i<j<k}j_{ijk*} \Omega^2_{X_{ijk}} \xrightarrow{d_2} \ldots
\end{gather}
Then a cocyle $c\in H^0(X,\prod_{i<j<k}j_{ijk*} \Omega^2_{X_{ijk}})$
defines a gerbe, with objects $a$ on $U \subset X$ given by
$a\in d_1^{-1}(c|_U)$ and morphisms 
$$Hom(a,b)= d_0^{-1}(b-a).
$$ 
We consider the image gerbe $c'\in
H^0(X,\prod_{i<j<k}j_{ijk*} (f^*\Omega^1_S \otimes
\Omega^1_{X/S})|_{X_{ijk}})$. In particular, for $X_i$ the
complement of \'etale multisections on a sufficiently small affine open
sets $S_{\alpha} \subset S$, the residue ($k_I:S_I \inj S$)
$${\rm res}(c') \in H^0(S,\prod_{\alpha < \beta } k_{\alpha\beta*}
\Omega^1_{S_{\alpha\beta}})
$$ 
represents the image torsor
\begin{gather}
0 \to \Omega^1_S \to d_0^{-1}({\rm res}(c')) \in 
\prod_{\alpha}k_{\alpha*}\Omega^1_{S_\alpha} \to {\rm
res}(c')\to 0
\end{gather} 
where 
\begin{gather}
\Omega^1_S \to \prod_{\alpha} k_{\alpha*}\Omega^1_{S_\alpha}
\xrightarrow{d_0} \prod_{\alpha < \beta} k_{\alpha
\beta*}\Omega^1_{S_{\alpha\beta}} \xrightarrow{d_1} \ldots
\end{gather}
is the Cech resolution of $\Omega^1_S$.

Assume now that $c'$ is a coboundary, and choose
\begin{gather}
\gamma\in H^0(X,\prod_{i<j}j_{ij*}(f^*\Omega^1_S\otimes
\Omega^1_{X/S})|_{X_{ij}}) 
\end{gather} 
with
$d_1(\gamma)=c$. So $\gamma$ is a trivialization of the gerbe $c'$
and ${\rm res} (\gamma) \in d_0^{-1}({\rm res}(c'))$ is a
trivialization of the torsor ${\rm res}(c')$. 
In particular, for $c$ equal, as above, $c_1^2-2c_2$ or $c_1^2$,
the torsor ${\rm res}(c')$ is the Atiyah torsor of the line
bundle  $f_*c \in \text{Pic}(S)$. Thus the choice of $\gamma$ defines
a splitting ${\rm res} (\gamma)$ of the Atiyah torsor, or
equivelently, it defines a
connection on the line bundle $f_*(c)$.

It remains to see why the choice of a relative connection
$\nabla_{X/S}: E \to \Omega^1_{X/S} \otimes E$ defines a
$\gamma$. The gauge transformation equation reads
\begin{gather}
dg_{ij}g^{-1}_{ij}= A_i - g_{ij}A_jg_{ij}^{-1} + \rho_{ij}
\end{gather}
where $A_i \in H^0(X_i, \Omega^1_X \otimes {\sE}nd(E))$
is a lifting to global forms of the local relative forms of the
connection, and $\rho_{ij} \in H^0(X_{ij},f^*(\Omega^1_S)\otimes
\Omega^1_{X/S}\otimes \sE nd(E))$ is a cocyle. 
This defines in natural way 
\begin{gather}
{\rm Tr} (dg_{ij}g_{ij}^{-1} dg_{jk}g_{jk}^{-1}) =
\delta(\eta_{ij})\\ 
{\rm Tr} (dg_{ij}g_{ij}^{-1}){\rm Tr}(dg_{jk}g_{jk}^{-1})=
\delta (\xi_{ij})
\end{gather} with
\begin{multline}\label{5.14}
\eta_{ij}= {\rm Tr}(A\cup \delta A+[A,\rho])_{ij} \\
={\rm
Tr}(A_i\delta(A)_{ij} -\rho_{ij}A_j + A_i\rho_{ij} )
\in H^0(X_{ij}, f^*\Omega^1_S \otimes \Omega^1_{X/S})\\
\xi_{ij}= (a_i\delta(a)_{ij} -\kappa_{ij}a_j + a_i\kappa_{ij}) 
\in H^0(X_{ij}, f^*\Omega^1_S \otimes \Omega^1_{X/S})\hfil
 \end{multline}
for $a= {\rm Tr}(A), \kappa= {\rm Tr}(\rho)$. 

We  claim  now that this construction leads to the same class as
ours. It is obvious for $f_*(c_1^2)$ because of the product
defined in section 3. As for $f_*(c_1^2-2c_2)$, we can exhibit the $AD$
class as follows. The relative $AD$-complex in this case is
$$\sK_2 \to \Omega^2_X/f^*\Omega^2_S.
$$
To write down a cocycle, we fix a flag in $E$, so the transition matrices
$g_{ij}$ are upper triangular:
$$g_{ij} = \begin{pmatrix}\ell^{(1)}_{ij} & \ldots \\
0 &\ell^{(2)}_{ij} & \ldots \\
\vdots & \vdots \\
0 & 0 & \ldots\ell^{(n)}_{ij}\end{pmatrix},
$$
and we take $a_{ijk} := \sum_p \{\ell^{(p)}_{ij},\ell^{(p)}_{jk}\}\in
\sK_2(U_{ijk})$. Clearly,
$$d\log(a_{ijk}) = \text{Tr}(dg_{ij}g_{ij}^{-1} dg_{jk}g_{jk}^{-1}).
$$
The claim is now that 
the $AD$-class is  represented by the $2$-hypercocycle
\begin{gather} \label{cocclaim}
(a_{ijk},\eta_{ij}) \in \Big(\sC^2(\sK_2) \times
\sC^1(f^*\Omega^1_S\otimes \Omega^1_{X/S})\Big)
\end{gather}
 where $\eta_{ij}$ is the same cochain \eqref{5.14}
which trivializes the gerbe $c'$ in Beilinson's construction. 

To see this, let us introduce the complete flag variety
$q: Q\to X$, such that $q^*E$ carries a tautological complete
flag. The construction of $\Omega$ (see \ref{Omega}) carries over
to the complete flag variety, and the connection $q^*\nabla$
on $q^*E$ induces  
a $\Omega^1_Q/\Omega$ connection, which then stabilizes the flag.
The cohomology group in which the class considered is living is
$AD^2_\tau(Q)$, and since $\Omega^3_Q/\Omega\cdot
\Omega^2_Q=q^*\Omega^3_{X/S} =0$ for a family of curves, one
has $AD^2_\tau(Q)= \H^2(Q, \sK_2 \xrightarrow{d\log }
\Omega^2_Q/\wedge^2 \Omega)$.

The equation of the connection becomes
\begin{gather}
dh_{ij}h_{ij}^{-1}= B_i-h_{ij}B_jh_{ij}^{-1} + \sigma_{ij}
\end{gather}
where $h_{ij}$ is an upper triangular matrix of functions,
$B_i$ is an upper triangular matrix of forms in $\Omega^1_Q$, 
$\sigma_{ij}$ is an upper triangular matrix of forms in $\Omega$. 
In particular, since $h$, $B$ and $\sigma$ are upper triangular,
the rule for the product implies that $q^*(c_1^2-2c_2)$ is
\begin{gather} \label{cocclaim3}
(b_{ijk}, \text{Tr}(B\cup \delta(B) + [B, \sigma])\in
\Big(\sC^2(\sK_2)\times \sC^1(\Omega^2_Q/\wedge ^2 \Omega)\Big) ,
\end{gather} 
where $b_{ijk}= \sum_p \{ \lambda_{ij}^{(p)},
\lambda_{jk}^{(p)} \}$ and where the $\lambda_{ij}^{(p)}$ are the
invertible diagonal entries of $h_{ij}$.

To finish the argument, one just needs that  \eqref{cocclaim3}
is compatible with gauge transformation.
\begin{remark}
The construction of a connection on $f_*c_2(E)$ is used in
\cite{F},
Theorem IV.3 in Faltings' construction of Hitchin's connection,
when $f$ is the projection $C \times M_\nabla \to M_\nabla$
to the moduli (or stack) of connections on a fixed curve $C$, and 
$E$ is the universal bundle. In particular, it identifies
$M_\nabla$ with the Atiyah torsor of $f_*c_2(E_0)$, where
$f_0: C \times M \to M$ is the projection to the moduli of
stable bundles (see Lemma IV.4 of \cite{F}), and $E_0$ is the
universal bundle on $C\times M$.
\end{remark}

\subsection{The trace complex} 

Let $f: X \to S$ be a smooth family of curves, and
let $E$ be a vector bundle on $X$ with a relative connection
$\nabla_{X/S}$. Under what conditions can one construct a connection on
$\det R f_*(E)$? Intuitively, the Riemann Roch theorem gives $[\det R
f_*(E)]\in
\text{Pic}(S)$ as a linear combination of 
$$f_*c_2(E),\ f_*(c_1(E)^2),
f_*(c_1(E)\cdot K_{X/S}), \text{ and } f_*(K_{X/S}^2). 
$$ 
We have shown
how to put connections on $f_*c_2(E)$ and $f_*(c_1(E)^2)$. Suppose in
addition that 
$X=Y\times S \to S$ is a product family so $K_{X/S}^2$ is trivial, and
assume further that we have an absolute connection $D$ on $\det E$ (e.g.
$\det E = \sO_X$). Intuitively again, $D$ gives a connection on
$f_*(c_1(E)\cdot K_{X/S})$. We sketch how in this situation one may use
the trace complex of Beilinson and Schechtman \cite{BS} (see also
\cite{ET}) to define a connection on $\det R f_*(E)$. 

Associated to the family $X/S$ and the bundle $E$ one has the Atiyah
algebra 
\eq{5.18}{ 0 \to \sE nd(E) \to \sA_{E,X} \to T_X \to 0.
}
One has a filtration $T_{X/S} \subset T_f \subset T_X$, where $T_f$
consists of vector fields whose image falls in $f^{-1}T_S \subset
f^*T_S$. By pullback one defines subalgebras 
\eq{5.19}{\sE nd(E) \subset \sA_{E,X/S} \subset \sA_{E,f} \subset
\sA_{E,X}.  } 
Pushing out by the trace map yields that corresponding Atiyah algebra for
$\det E$
\eq{5.20}{\begin{CD} 0 @>>> \sE nd(E) @>>> \sA_{E,?} @>>> T_? @>>> 0 \\
@. @VV \Tr V @VVV @| \\
0 @>>> \sO_X @>>> \sA_{\det E,?} @>>> T_? @>>> 0.
\end{CD}
} 
Thus we get
\eq{5.21}{0 \to \sE nd^0(E) \to \sA_{E,?}\to \sA_{\det E,?}\to 0. 
}

The trace complex
\eq{5.22}{{}^{tr}\!\!\sA^\bullet = \{\sA^{-2} \to \sA^{-1} \to \sA^0 \} 
}
fits into an exact sequence of complexes
\eq{5.23}{ 0 \to \omega_{X/S}[2] \to {}^{tr}\!\!\sA^\bullet \to
\{\sA_{E,X/S}
\to \sA_{E,f}\} \to 0 }
One has
\eq{5.24}{\sA^{-2} = \sO_X;\quad \sA^0 = \sA_{E,X/S}
}
and $\sA^{-1}$ is an extension
\eq{5.25}{\begin{CD} 
0 @>>> \omega_{X/S} @>>> \sA^{-1} @>>> \sA_{E,X/S} @>>> 0. 
\end{CD}
}

An absolute connection on $\det E$ induces compatible splittings
\eq{5.26}{\begin{array}{ccc}\sA_{\det E,X/S} &
\stackrel{\leftarrow}{\twoheadrightarrow}& T_{X/S} \\
\downarrow && \downarrow \\\sA_{\det E,f} &
\stackrel{\leftarrow}{\twoheadrightarrow}& T_{f}.
\end{array}
}
Assuming $X=Y\times S$ (or, more generally, that we are given a
connection for the map $f:X \to S$) we get a decomposition
\eq{5.27}{ T_f = T_{X/S}\oplus f^{-1}T_S
}
Combining \eqref{5.26} and \eqref{5.27} yields an injective
quasiisomorphism of complexes
\eq{5.28}{ \{0 \to f^{-1}T_S\} \inj \{\sA_{\det E, X/S} \to \sA_{\det E,
f}\}. }
We can pull back the sequences \eqref{5.21} along this map to get a
quasiisomorphic subcomplex (defining $\sB_{\det E,f}$)
\eq{5.29}{\{\sE nd(E)^0 \to \sB_{\det E,f} \} \inj\{\sA_{E,X/S} \to
\sA_{E,f}\} }
Finally, we can pull back \eqref{5.23} along this map, defining a
quasiisomorphic subcomplex ${}^{tr}\sB^\bullet\subset
{}^{tr}\!\!\sA^\bullet$ with
\begin{gather}\sB^{-2}=\sA^{-2}=\sO_X;\quad \sB^{0} = \sB_{\det E,f} \\
\label{5.31} 0 \to \omega_{X/S} \to \sB^{-1} \to \sE nd(E)^0 \to 0
\end{gather}
The sequence \eqref{5.31} is the trace-free part of the dual of the
relative Atiyah sequence, so it is split by the given relative connection
on $E$. there results a map of complexes
\eq{5.32}{{}^{tr}\sB^\bullet \to \omega^\bullet[2]
}
splitting the sequence \eqref{5.23} upto quasiisomorphism. Applying
$R^0f_*$ yields a splitting for the exact sequence
\eq{}{0 \to \sO_S\cong R^2f_*\omega^\bullet \to R^0f_*{}^{tr}\sB^\bullet
\to T_S \to 0.  
}
But, by \cite{BS}, this is the Atiyah sequence for $\det Rf_*(E)$, and a
splitting defines a connection on that line bundle. In summary, we have
shown
\begin{thm} Let $f:X\to S$ be a trivial family of curves. Let
$E$ be a bundle on $X$, together with a global connection 
$D: {\rm det}(E) \to \Omega^1_X \otimes{\rm det}(E)$ and a
relative one
$\nabla_{X/S}: E\to \Omega^1_{X/S} \otimes E.$
Then $D$ and  $\nabla_{X/S}$ define a subcomplex ${}^{tr}\sB^{\bullet}
\subset {}^{tr}\!\!\sA^{\bullet}$ of the trace complex,
quasiisomorphic to it, and a splitting
${}^{tr}\sB^{\bullet} \to \omega[2]$. In particular, these data induce a
connection on the determinant of $Rf_*(E)$. 
\end{thm}

We admit to not having worked out the precise relation between this
connection and the connections on $f_*(c_2)$ and $f_*(c_1^2)$ described
earlier.

\newpage
\bibliographystyle{plain}
\renewcommand\refname{References}

\end{document}